\documentclass{amsart}
 \usepackage{amssymb}

\DeclareMathSymbol{\twoheadrightarrow} {\mathrel}{AMSa}{"10}

\def\Q{{\mathbf Q}}
\def\Z{{\mathbf Z}}

\def\F{{\mathbf F}}

\def\Sn{{\mathbf S}_n}
\def\An{{\mathbf A}_n}

\def\Gal{\mathrm{Gal}}

\def\Perm{\mathrm{Perm}}

\def\pr{\mathrm{pr}}

\def\pr{\mathrm{pr}}

\def\End{\mathrm{End}}

\def\Aut{\mathrm{Aut}}
\def\Hom{\mathrm{Hom}}

    \def\RR{\mathfrak{R}}

\def\fchar{\mathrm{char}}

\def\GL{\mathrm{GL}}

        \def\res{\mathrm{res}}

\def\dim{\mathrm{dim}}
\def\O{{\mathfrak O}}

\newtheorem{thm}{Theorem}[section]

\newtheorem{lem}[thm]{Lemma}

\newtheorem{cor}[thm]{Corollary}

 \newtheorem{defn}[thm]{Definition}

\newtheorem{rem}[thm]{Remark}

\begin{document}

\title{Non-isogenous superelliptic Jacobians}
\author{Yuri G. Zarhin}
\address{Department of Mathematics, Pennsylvania
State University, University Park, PA 16802, USA}
\email{zarhin\char`\@math.psu.edu}

\maketitle

\begin{abstract}
Let $\ell$ be an odd prime. Let $K$ be a field of characteristic
zero with algebraic closure $K_a$. Let $n,m \ge 4$ be integers
that are not divisible by $\ell$. Let $f(x), h(x) \in K[x]$ be
irreducible separable polynomials of degree $n$ and $m$
respectively. Suppose that the Galois group $\Gal(f)$ of $f$ acts
doubly transitively on the set $\RR_f$ of roots of $f$ and that
$\Gal(h)$  acts doubly transitively on  $\RR_h$ as well.
  Let $J(C_{f,\ell})$ and
$J(C_{h,\ell})$ be the Jacobians of the superelliptic curves
$C_{f,\ell}:y^{\ell}=f(x)$ and $C_{h,\ell}:y^{\ell}=h(x)$
respectively.
 We prove that $J(C_{f,\ell})$  and $J(C_{h,\ell})$ are not
isogenous over $K_a$ if the splitting fields of $f$ and $h$ are
linearly disjoint over $K(\zeta_{\ell})$.
\end{abstract}

\section{Definitions, notations, statements}

Let $K$ be a field. Let us fix its algebraic closure $K_a$ and
denote by $\Gal(K)$ the absolute Galois group  $\Aut(K_a/K)$ of
$K$. If $L \supset K$ is an overfield of $K$ and $L_a$ contains
$K_a$ (i.e., $K_a$ is the algebraic closure of $K$ in $L_a$) then
$K_a$ is  $\Aut(L_a/L)$-stable and we write
$$\res(L,K):\Gal(L)=\Aut(L_a/L)\to \Aut(K_a/K)=\Gal(K)$$
for the corresponding restriction map. If $X$ is an abelian
variety over $K_a$ then we write $\End(X)$ for the ring of all its
$K_a$-endomorphisms;  $1_X$ stands for the identity automorphism
of $X$. If $Y$ is an abelian variety over $K_a$ then we write
$\Hom(X,Y)$ for the (free commutative) group of all
$K_a$-homomorphisms from $X$ to $Y$. It is well-known that
$\Hom(X,Y)=0$ if and only if $\Hom(Y,X)=0$. If $X$ is defined over
$K$ then $X(K_a)$ carries a natural structure of $\Gal(K)$-module.
One may also view $X$ as an abelian variety over $L$; the subgroup
$X(K_a)\subset X(L_a)$ is $\Gal(L)$-stable and the corresponding
homomorphism $\Gal(L)\to \Aut(X(K_a))$ is the composition of
$\res(L,K):\Gal(L)\to \Gal(K)$ and the structure homomorphism
$\Gal(K)\to\Aut(X(K_a))$.

Let $f(x) \in K[x]$ be a polynomial of degree $n\ge 4$ without
multiple roots.  We write $\RR_f\subset K_a$ for the set of its
roots,  $K(\RR_f) \subset K_a$ for the splitting field of $f$ and
$\Gal(f)=\Aut(K(\RR_f)/K)=\Gal(K(\RR_f)/K)$ for the Galois group
of $f$. Then $\RR_f$ consists of $n=\deg(f)$ elements. The group
$\Gal(f)$ permutes elements of $\RR_f$ and therefore can be
identified with a certain subgroup of the group $\Perm(\RR_f)$ of
all permutations of  $\RR_f$. Clearly, every ordering of $\RR_f$
provides an isomorphism between $\Perm(\RR_f)$ and the full
symmetric group $\Sn$ which makes $\Gal(f)$ a certain subgroup of
$\Sn$. (This permutation subgroup is transitive if and only if $f$
is irreducible over $K$.)

Let $\ell$ be an odd prime. We write $\Z[\zeta_{\ell}]$ for the
ring of all integers in the $\ell$th cyclotomic ring
$\Q[\zeta_{\ell}]$.

 Let us assume that $\fchar(K)\ne \ell$ and consider the
superelliptic curve
$$C_{f,\ell}: y^{\ell}=f(x),$$
defined over $K$. Its genus $g=g(C_{f,\ell})$ equals
$(n-1)(\ell-1)/2$ if $\ell$ does not divide $n$  and
$(n-2)(\ell-1)/2$ if $\ell\mid n$. Let $J(C_{f,\ell})$ be the
Jacobian of $C_f$; it is a $g$-dimensional abelian variety over
$K_a$ that is defined over $K$. Then $\End(J(C_{f,\ell}))$
contains a certain subring isomorphic to $\Z[\zeta_{\ell}]$ (see
Sect. 3.2).

The main result of the present paper is the following statement.

\begin{thm}[Main Theorem]
\label{main} Suppose that $K$ is a field of characteristic
different from $\ell$ that contains a primitive $\ell$th root of
unity.  Let $f(x), h(x)\in K[x]$ be separable irreducible
polynomials of degree $n\ge 4$ and $m\ge 4$ respectively. Suppose
that the splitting fields of
 $f$ and $h$ are linearly disjoint over  $K$.

Suppose that the following conditions hold:

\begin{enumerate}
\item[(i)] The group
$\Gal(f)$   acts doubly transitively on  $\RR_f$; if $\ell$
divides $n$ then this action is $3$-transitive.
\item[(ii)] The group
$\Gal(h)$   acts doubly transitively on  $\RR_h$; if $\ell$
divides $n$ then this action is $3$-transitive.
\end{enumerate}

Then either $$\Hom(J(C_{f,\ell}), J(C_{h,\ell}))=0,\quad
\Hom(J(C_{h,\ell}),J(C_{f,\ell}))=0$$ or $p:=\fchar(K)>0$ and
there exists an abelian variety $Z$ defined over an algebraic
closure $\bar{\F}_{p}$ of $\F_{p}$ such that both $J(C_{f,\ell})$
and $J(C_{h,\ell})$ are isogenous over $K_a$ to self-products of
$Z$.
\end{thm}

\begin{rem}
The case $\ell=2$ (of hyperelliptic Jacobians) was treated in
\cite{ZarhinSh,ZarhinL}. See \cite{Mortimer,DM} for the list of
known doubly transitive permutation groups.
\end{rem}

The paper is organized as follows. In Sections \ref {homav} and
\ref{disjoint} we study pairs of abelian varieties with
homomorphism groups of big rank. We prove Theorem \ref{main} in \S
\ref{proofm}. Sections \ref{rep} and \ref{proofA} contain the
proof of some auxiliary results.

I am grateful to the referee, whose comments helped to improve the
 exposition.

\section{Homomorphisms of abelian varieties: statements}
\label{homav}

First, we need to introduce some notions from the theory of
abelian varieties. Let $K$ be a field and $d$ be a positive
integer that is not divisible by $\fchar(K)$. Let $X$ be an
abelian variety
  of positive dimension
defined over $K$. We write $X_d$ for the kernel of multiplication
by $d$ in $X(K_a)$. The commutative group $X_d$ is a free
$\Z/d\Z$-module of rank $2\dim(X)$ \cite{Mumford}. Clearly, $X_d$
is a Galois submodule in $X(K_a)$ . We write
$$\tilde{\rho}_{d,X}:\Gal(K) \to \Aut_{\Z/d\Z}(X_d) \cong
\GL(2\dim(X),\Z/d\Z)$$ for the corresponding (continuous)
homomorphism defining the Galois action on $X_d$. Let us put
$$\tilde{G}_{d,X}=\tilde{\rho}_{d,X}(\Gal(K)) \subset
\Aut_{\Z/d\Z}(X_d).$$ Clearly,  $\tilde{G}_{d,X}$ coincides with
the Galois group of the field extension  $K(X_d)/K$ where $K(X_d)$
is the field of definition of all points of order dividing $d$ on
$X$. In particular, if
 $\ell\ne \fchar(K)$ is a prime then
$X_{\ell}$ is a $2\dim(X)$-dimensional vector space over the prime
field $\F_{\ell}=\Z/\ell\Z$ and the inclusion $\tilde{G}_{\ell,X}
\subset \Aut_{\F_{\ell}}(X_{\ell})$ defines a faithful linear
representation of the group  $\tilde{G}_{\ell,X}$ in the vector
space $X_{\ell}$.

We write $\End_K(X)\subset \End(X)$ for the (sub)ring of all
$K$-endomorphisms of $X$ and $\End_K^0(X)\subset \End^0(X)$ for
the corresponding $\Q$-(sub)algebra of all $K$-endomorphisms of
$X$. If $Y$ is  an abelian variety over $K$ then we write
$\Hom^0(X,Y)$ for the $\Q$-vector space $\Hom(X,Y)\otimes\Q$.

 Let $E$ be a number field and $\O \subset E$ be the ring of all
its  algebraic integers. Let $(X, i)$ be a pair consisting of an
abelian variety $X$ over $K_a$ and an embedding
$$i:E \hookrightarrow  \End^0(X)$$
such that $i(1)=1_X$. The degree $[E:\Q]$ divides $2\dim(X)$ (see
\cite{Ribet}).

If $r$ is a positive integer then we write $i^{(r)}$ for the
composition
$$E\hookrightarrow \End^0(X)\subset \End^0(X^r)$$
of $i$ and the diagonal inclusion $\End^0(X)\subset \End^0(X^r)$.

If $(Y,j)$ is a   pair consisting of an abelian variety $Y$ over
$K_a$ and an embedding $j:E \hookrightarrow  \End^0(Y)$ with
$j(1)=1_Y$ then we write
$$\Hom^0((X,i),(Y,j))=\{u\in \Hom^0(X,Y)\mid u i(e)=j(e) u \quad
\forall u \in E\}.$$ Clearly, $\Hom^0((X,i),(Y,j))$ carries a
natural structure of finite-dimensional $E$-vector space. Notice
that the $\Q$-vector space $\Hom^0(X,Y)$ carries a natural
structure of $E\otimes_{\Q}E$-module defined by the formula
$$(e_1\otimes e_2)\phi =j(e_1)\phi i(e_2) \quad \forall e_1,e_2\in
E, \phi\in \Hom^0(X,Y).$$

\begin{rem}
\label{normal} It is well-known that if the field extension
$E/\Q$ is normal then for each automorphism
$\sigma\in\Aut(E)=\Gal(E/Q)$ there is a surjective $E$-algebra
homomorphism
$$\pr_{\sigma}:E\otimes_{\Q}E\twoheadrightarrow E, e_1\otimes
e_2\mapsto e_1\sigma(e_2).$$ (Here the structure of $E$-algebra on
$E\otimes_{\Q}E$ is defined  by
$$e (e_1\otimes e_2)=e e_1\otimes e_2 \quad \forall e,e_1,e_2\in
E.$$ The well-known $E$-linear independence of all $\sigma:E \to
E$ implies that the direct sum of all $\pr_{\sigma}$'s is an
isomorphism
$$\oplus_{\sigma\in\Gal(E/\Q)}\pr_{\sigma}:E\otimes_{\Q}E=\oplus_{\sigma \in \Gal(E/\Q)}E.$$
This allows us to view $\pr_{\sigma}$ as mutually orthogonal
projection maps $\pr_{\sigma}:E\otimes_{\Q}E\to E\otimes_{\Q}E$,
whose sum is the identity map. Also, the annihilator of
$\sigma(e)\otimes 1-1\otimes e$ in $E\otimes_{\Q}E$ coincides with
the image $\pr_{\sigma}:E\otimes_{\Q}E$ of $\pr_{\sigma}$.

 This
implies easily that
$$\Hom^0((X,i),(Y,j\sigma))=\pr_{\sigma}(\Hom^0(X,Y))$$
and
$$\Hom^0(X,Y)=\oplus_{\sigma \in
\Gal(E/\Q)}\Hom^0((X,i\sigma),(Y,j))\eqno(1)$$ where $i\sigma :E
\hookrightarrow \End^0(X)$ is the composition of the automorphism
$\sigma:E\to E$ and $i:E \to\End^0(X)$.
\end{rem}

Let us denote by $\End^0(X,i)$ the centralizer of $i(E)$ in
$\End^0(X)$. Clearly, $\End^0(X,i)=\Hom^0((X,i),(X,i))$ and $i(E)$
lies in the center of the finite-dimensional $\Q$-algebra
$\End^0(X,i)$. It follows that $\End^0(X,i)$ carries a natural
structure of finite-dimensional $E$-algebra. One may easily check
\cite[Remark 4.1]{ZarhinL} that
 $\End^0(X,i)$ is a semisimple $E$-algebra; it is simple if and only if
$X$ is isogenous to a self-product of an (absolutely) simple
abelian variety. The following two
 assertions are contained in \cite[Theorem 4.2]{ZarhinL},
 \cite[Remark 3.2]{ZarhinM}.

\begin{thm}
\label{boundE} $\dim_E(\End^0((X,i)) \le
\frac{4\cdot\dim(X)^2}{[E:\Q]^2}$.
\end{thm}

\begin{thm}
\label{maxE} Suppose that
 $$\dim_E(\End^0((X,i)) =
\frac{4\cdot\dim(X)^2}{[E:\Q]^2}.$$
 Then:
\begin{itemize}
\item[(i)] $\End^0((X,i)$ is a central simple $E$-algebra.
\item[(ii)] There exists an absolutely simple abelian variety $B$
of CM-type over $K_a$ such that $X$ is isogenous to a self-product
of $B$. \item[(iii)] If $\fchar(K)=0$ then  $[E:\Q]$ is even and
there exist a $[E:\Q]/2$-dimensional abelian variety $Z$ over
$K_a$, an isogeny $\psi: Z^r \to X$ and an embedding $k: E
\hookrightarrow \End^0(Z)$ that send $1$ to $1_Z$ and such that
$\psi \in \Hom^0((Z^r,k^{(r)}),(X,i))$.
\end{itemize}
\end{thm}

\begin{rem}
\label{CM} Suppose that
 $$\dim_E(\End^0((X,i)) =
\frac{4\cdot\dim(X)^2}{[E:\Q]^2}$$ and $\fchar(K)>0$. In notations
of Theorem \ref{maxE}, it follows from a a theorem of Grothendieck
\cite[Th. 1.1]{Oort} that $B$ is isogenous to an abelian variety
defined over a finite field. This implies that $X$ is also
isogenous to an abelian variety defined over a finite field.
\end{rem}

If $i(\O) \subset \End(X)$ and $j(\O) \subset \End(Y)$ then we put
$$\Hom((X,i),(Y,j))=\{u\in \Hom(X,Y)\mid u i(e)=j(e) u \quad
\forall u \in E\}.$$ Clearly,
$$\Hom^0((X,i),(Y,j))=\Hom((X,i),(Y,j))\otimes \Q,$$
$$ \Hom((X,i),(Y,j))=\Hom^0((X,i),(Y,j))\bigcap \Hom(X,Y),$$
which is   an $\O$-lattice in the $E$-vector space
$\Hom^0((X,i),(Y,j))$.

\begin{rem}
\label{power} There are canonical isomorphisms of $E$-vector
spaces
$$\Hom^0((X^r,i^{(r)}),(Y,j))=(\Hom^0((X,i),(Y,j)))^r=\Hom^0((X^r,i),(Y,j^{(r)}))$$
where $(\Hom^0((X,i),(Y,j)))^r$ is a direct sum of $r$ copies of $\Hom^0((X,i),(Y,j))$.
It follows easily that there is a canonical isomorphism of $E$-vector spaces
$$\Hom^0((X^r,i^{(r)}),(Y,j^{(m)}))=(\Hom^0((X,i),(Y,j)))^{rm}$$
for all positive integers $r$ and $m$.
\end{rem}

\begin{lem}
\label{dimb}
\begin{itemize}
\item[(i)] $\dim_E(\Hom^0((X,i),(Y,j))) \le
\frac{4\cdot\dim(X)\dim(Y)}{[E:\Q]^2}$;
 \item[(ii)] If
$\dim(X)=\dim(Y)$ and
$$\dim_E(\Hom^0((X,i),(Y,j))) =
\frac{4\cdot\dim(X)\dim(Y)}{[E:\Q]^2}$$ then
 $\Hom^0((X,i),(Y,j))$ contains an isogeny $\phi:X \to Y$. In particular,
$$\Hom^0((X,i),(Y,j))=\phi\cdot \End^0(X,i), \
\End^0(Y,j)=\phi\End^0(X,i)\phi^{-1}$$ and
$$\dim_E\End^0(Y,j)=\dim_E\End^0(X,i)
=\dim_E\Hom^0((X,i),(Y,j))=$$
$$\frac{4\dim(X)^2}{[E:\Q]^2}=\frac{4\dim(Y)^2}{[E:\Q]^2}.$$
\end{itemize}
\end{lem}

\begin{proof}[of Lemma \ref{dimb}]
Let us fix a prime $\ell \ne \fchar(K)$. Let us put
$$E_{\ell}:=E\otimes_{\Q}\Q_{\ell}.$$
Clearly, $E_{\ell}$ is a direct sum of finitely many $\ell$-adic
fields.

 Let $T_{\ell}(X)$ be the $\Z_{\ell}$-Tate module of $X$
\cite{Mumford}. Recall that $T_{\ell}(X)$ is a free
$\Z_{\ell}$-module of rank $2\dim(X)$. Let us put
$$V_{\ell}(X)=T_{\ell}(X)\otimes_{\Z_{\ell}}\Q_{\ell};$$
it is a $2\dim(X)$-dimensional $\Q_{\ell}$-vector space. There are
natural embeddings
$$\End^0(X)\otimes_{\Q}\Q_{\ell}\hookrightarrow
\End_{\Q_{\ell}}V_{\ell}(X), \
\End^0(Y)\otimes_{\Q}\Q_{\ell}\hookrightarrow
\End_{\Q_{\ell}}V_{\ell}(Y),$$
$$\Hom^0(X,Y)\otimes_{\Q}\Q_{\ell}\hookrightarrow
\Hom_{\Q_{\ell}}(V_{\ell}(X),V_{\ell}(Y)).$$ Now the injections
$i$ and $j$ give rise to the injections
$$E_{\ell}\hookrightarrow \End^0(X)\otimes_{\Q}\Q_{\ell}\hookrightarrow
\End_{\Q_{\ell}}V_{\ell}(X), \ E_{\ell}\hookrightarrow
\End^0(Y)\otimes_{\Q}\Q_{\ell}\hookrightarrow
\End_{\Q_{\ell}}V_{\ell}(Y).$$ These injections provide
 $V_{\ell}(X)$ and $V_{\ell}(Y)$ with the natural structure of free
$E_{\ell}$-modules of rank $\frac{2\dim(X)}{[E:\Q]}$ and
$\frac{2\dim(Y)}{[E:\Q]}$ respectively \cite{Ribet}. Clearly, the
image of
$$\Hom^0((X,i),(Y,j)\otimes_{\Q}\Q_{\ell}\subset
\Hom^0(X,Y)\otimes_{\Q}\Q_{\ell}$$ in
$\Hom{\Q_{\ell}}(V_{\ell}(X),V_{\ell}(Y))$ lies in
$\Hom_{E_{\ell}}(V_{\ell}(X),V_{\ell}(Y))$; in fact, it is a free
$E_{\ell}$-submodule of $\Hom_{E_{\ell}}(V_{\ell}(X),V_{\ell}(Y))$
of rank  $\dim_E(\Hom^0((X,i),(Y,j))$. The rank of the free
$E_{\ell}$-module $\Hom_{E_{\ell}}(V_{\ell}(X),V_{\ell}(Y))$
equals  the product of the ranks of $V_{\ell}(X)$ and
$V_{\ell}(Y)$, i.e. equals
$$\frac{2\dim(X)}{[E:\Q]} \cdot \frac{2\dim(Y)}{[E:\Q]}.$$
We conclude that
$$\dim_E(\Hom^0((X,i),(Y,j))\le \frac{2\dim(X)}{[E:\Q]} \cdot
\frac{2\dim(Y)}{[E:\Q]}.$$ Clearly, the equality holds if and only
if
$$\Hom_{E_{\ell}}(V_{\ell}(X),V_{\ell}(Y))=\Hom^0((X,i),(Y,j))\otimes_{\Q}\Q_{\ell}.$$
Suppose that the equality holds and assume, in addition, that
$\dim(X)=\dim(Y)$. Then the ranks of $V_{\ell}(X)$ and
$V_{\ell}(Y)$ do coincide and there exists an isomorphism $u:
V_{\ell}(X)\cong V_{\ell}(Y)$ of $E_{\ell}$-modules. Since $\Q$ is
everywhere dense in $\Q_{\ell}$ in the $\ell$-adic topology, there
exists $u'\in \Hom^0((X,i),(Y,j))$ that is also an isomorphism
between $V_{\ell}(X)$ and $V_{\ell}(Y)$. Replacing $u'$ by $N u'$
for suitable positive integer $N$, we may assume that $u' \in
\Hom(X,Y)$. Then $u'$ must be an isogeny.
\end{proof}

\begin{thm}
\label{manyhom} Suppose that $E$ is a number field, $X$ and $Y$
are  abelian varieties of positive dimension over an algebraically
closed field $K_a$,
$$i:E \hookrightarrow \End^0(X), \quad j:E \hookrightarrow \End^0(Y)$$
are embeddings that send $1$ to the identity automorphisms of $X$
and $Y$ respectively. Let us put
$$r_X:=\frac{2\dim(X)}{[E:\Q]}, \quad r_Y:=\frac{2\dim(Y)}{[E:\Q]}.$$
Let us assume that
$$\dim_E \Hom^0((X,i),(Y,j))=r_X\cdot r_Y.$$
Then both $\End^0(X,i)$ and $\End^0(Y,j)$ are central simple
$E$-algebras and
$$\dim_E\End^0(X,i)=r_X^2,
\quad \dim_E\End^0(Y,j)=r_Y^2.$$ In addition,
 both $X$ and $Y$ are  isogenous to self-products of a
certain absolutely simple abelian variety $B$ of CM-type.
\end{thm}

\begin{proof}[of Theorem \ref{manyhom}]
Clearly,
$$\dim(X^{r_Y})=\frac{2\dim(X)\dim(Y)}{[E:\Q]}=\dim(Y^{r_X}).$$
It follows from Remark \ref{power} that
$$\dim_E (\Hom^0((X^{r_Y},i^{(r_Y)}),(Y^{r_X},j^{(r_X)}))=
\frac{4\dim(X^{r_Y})\dim(Y^{r_X})}{[E:\Q]}.$$

By Lemma \ref{dimb}  there exists an isogeny $\phi:X^{r_Y}\to
Y^{r_X}$ that lies in
$$\Hom^0((X^{r_Y},i^{(r_Y)}),(Y^{r_X},j^{(r_X)}).$$
In addition,  $\End^0(X,i)$ and $\End^0(Y,j)$ are central simple
$E$-algebras and
$$\dim_E\End^0(X^{r_Y},i^{(r_Y)})=\left(\frac{2\dim(X^{r_Y})}{[E:\Q]}\right)^2
=(r_X r_Y)^2.$$ Similarly,
$$\dim_E\End^0(Y^{r_X},j^{(r_X)})=(r_Y r_X)^2.$$
This implies the first claim.

 Applying Theorem \ref{maxE} to both
$(X,i)$ and $(Y,j))$, we conclude that there exist absolutely
simple abelian varieties $B$ and say, $B'$ of CM-type such that
$X$ is  isogenous to a self-product of $B$ and $Y$ is  isogenous
to a self-product of $B'$.
 Since $\Hom^0((X,i), (Y,j))\ne 0$, we conclude that $\Hom(X,Y)\ne 0$ and therefore $\Hom(B,B')\ne
 0$. This implies that $B$ and $B'$ are isogenous and therefore $Y$ is  isogenous
to a self-product of $B$.
\end{proof}

Suppose that $X$ is defined over $K$ and $i(\O) \subset
\End_K(X)$. Then we may view elements of $\O$ as $K$-endomorphisms
of $X$.

 Let $\lambda$ be a maximal ideal in $\O$. We write $k(\lambda)$
for the corresponding (finite) residue field. Let us put
$$X_{\lambda}=X_{\lambda,i}:=\{x \in X(K_a)\mid i(e)x=0 \quad \forall e\in \lambda\}.$$
Clearly, if $\fchar((k)(\lambda))=\ell$  then $\lambda\supset \ell
\cdot\O$ and therefore $X_{\lambda}\subset X_{\ell}$. Moreover,
$X_{\lambda}$ is a Galois submodule of $X_{\ell}$ and
$X_{\lambda}$ carries a natural structure of
$\O/\lambda=k(\lambda)$-vector space. It is known \cite{Ribet}
that if $\ell \ne \fchar(K)$ then
$$\dim_{k(\lambda)}X_{\lambda}= \frac{2\dim(X)}{[E:\Q]}.$$
We write
$$\tilde{\rho}_{\lambda,X}=\tilde{\rho}_{\lambda,X,K}:\Gal(K) \to \Aut_{k(\lambda)}(X_{\lambda}) \cong
\GL(d_{X,E},k(\lambda))$$ for the corresponding (continuous)
homomorphism defining the Galois action on $X_{\lambda}$.
 Let us
put
$$\tilde{G}_{\lambda,X}=\tilde{G}_{\lambda,i,X}:=\tilde{\rho}_{\lambda,X}(\Gal(K)) \subset
\Aut_{k(\lambda)}(X_{\lambda}).$$ Clearly, $\tilde{G}_{\lambda,X}$
coincides with the Galois group of the field extension
$K(X_{\lambda})/K$ where $K(X_{\lambda})=K(X_{\lambda,i})$ is the
field of definition of all points in $X_{\lambda}$.

In order to describe $\tilde{\rho}_{\lambda,X,K}$ explicitly, let
us assume for the sake of simplicity that $\lambda$ is the only
maximal ideal of $\O$ dividing $\ell$, i.e.,
$\ell\cdot\O=\lambda^b$ where the positive integer $b$ satisfies
$[E:\Q]=b \cdot \dim_{\F_{\ell}}k(\lambda)$. Then
$\O\otimes\Z_{\ell}=\O_{\lambda}$ where $\O_{\lambda}$ is the
completion of $\O$ with respect to $\lambda$-adic topology. Let us
choose an element $c \in \lambda$ that does not lie in
$\lambda^2$. One may easily check \cite[\S 3]{ZarhinM} that
$$X_{\lambda}=\{x\in X_{\ell}\mid cx=0\}\subset X_{\ell}.$$

Let $T_{\ell}(X)$ be the $\Z_{\ell}$-Tate module of $X$. Recall
that $T_{\ell}(X)$ is a free $\Z_{\ell}$-module of rank $2\dim(X)$
provided with the continuous action
$$\rho_{\ell,X}:\Gal(K) \to \Aut_{\Z_{\ell}}(T_{\ell}(X))$$
and the natural embedding
$$\End_K(X)\otimes\Z_{\ell} \hookrightarrow \End_{\Z_{\ell}}(T_{\ell}(X)),$$
whose image commutes with $\rho_{\ell,X}(\Gal(K))$. In particular,
$T_{\ell}(X)$ carries the natural structure of
$\O\otimes\Z_{\ell}=\O_{\lambda}$-module; it is known \cite{Ribet}
that the $\O_{\lambda}$-module $T_{\ell}(X)$ is free of rank
$d_{X,E}$. There is also the natural isomorphism of Galois modules
$$X_{\ell}=T_{\ell}(X)/\ell T_{\ell}(X),$$
which is also an isomorphism of $\End_K(X)\supset \O$-modules. One
may easily check \cite[\S 3]{ZarhinM} that  the
$\O[\Gal(K)]$-module
$$X_{\lambda}=T_{\ell}(X)/(\lambda\O_{\lambda})T_{\ell}(X)=T_{\ell}(X)\otimes_{\O_{\lambda}}k(\lambda).$$

\begin{rem}
\label{aut} Let $\sigma$ be an automorphism of $E$. Clearly,
$\sigma(\O)=\O$ and $\sigma(\lambda)=\lambda$ (since $\lambda$ is
the only maximal ideal dividing $\ell$). However, $\sigma$ may
induce a non-trivial automorphism of $k(\lambda)$ (if $k(\lambda)
\ne \F_{\ell}$). Let us consider the composition
$$t:=i\sigma: E \hookrightarrow \End_K^0(X).$$
Clearly, $t(\O)=i(\O)\subset \End_K(X)$. It is also clear that
$$X_{\lambda}=X_{\lambda,i}=X_{\lambda,t}, \
K(X_{\lambda})=K(X_{\lambda,i})=K(X_{\lambda,t}), \
\tilde{G}_{\lambda,X}=\tilde{G}_{\lambda,i,X}=\tilde{G}_{\lambda,t,X}.$$
However, the structure of the $k(\lambda)$-vector space on
$X_{\lambda,t}$ is the twist via $\sigma$ of   the structure of
the $k(\lambda)$-vector space on $X_{\lambda,i}$. This means that
multiplication by any $a\in k(\lambda)$ in $X_{\lambda,t}$
coincides with multiplication by $\sigma(a)$ in $X_{\lambda,i}$.
However, this {\sl twist} does not change the algebra of linear
operators, i.e.
 $$\End_{k(\lambda)}(X_{\lambda,i})=\End_{k(\lambda)}(X_{\lambda,t}),
 \ \Aut_{k(\lambda)}(X_{\lambda,i})=\Aut_{k(\lambda)}(X_{\lambda,t}).$$
This implies that the centralizers of $\tilde{G}_{\lambda,X}$ in
 $\End_{k(\lambda)}(X_{\lambda,i})$ and
 $\End_{k(\lambda)}(X_{\lambda})$ do coincide.
 In particular, if the centralizer
 $\End_{\tilde{G}_{\lambda,X}}(X_{\lambda})$
 is $k(\lambda)$ (resp. a field) then the centralizer
 $\End_{\tilde{G}_{\lambda,X}}(X_{\lambda,i})$
 is also $k(\lambda)$ (resp. a field).
\end{rem}

\begin{rem}
\label{basechange}
 Suppose that $L$ is an overfield of $K$ and
$K_a$ is the algebraic closure of $K$ in $L_a$. Then one may view
$X$ as an abelian variety over $L$ and $i(\O)\subset \End_L(X)$.
The base change (from $K$ to $L$) does not change the groups $X_n$
and $X_{\lambda}$. One may easily check that
$\tilde{\rho}_{\lambda,X,L}:\Gal(L) \to
\Aut_{k(\lambda)}(X_{\lambda})$ coincides with the composition of
$\res(L,K):\Gal(L)\to\Gal(K)$ and
$\tilde{\rho}_{\lambda,X,K}:\Gal(K) \to
\Aut_{k(\lambda)}(X_{\lambda})$.
\end{rem}

\section{Disjoint abelian varieties}
\label{disjoint} Throughout this Section
 $E$ is a number field with the ring of integers $\O$ and  $\lambda$ is a maximal ideal in $\O$,
whose residue field $k(\lambda)=\O/\lambda$ has characteristic
$\ell$. We assume that $\lambda$ is the only maximal ideal of $\O$
dividing $\ell$. Let $K$ a field of characteristic different from
$\ell$. Let $X$ and $Y$ are
 abelian varieties of positive dimension
over $K$ provided with embeddings
$$i:E\hookrightarrow \End^0_K(X)\subset \End^0(X), \quad j:E\hookrightarrow
\End_K^0(Y)\subset \End^0(Y)$$ such that
 $$1_X=i(1)\in i(\O)\subset \End_K(X), \ 1_Y=j(1)\in j(\O)\subset \End_K(Y).$$
 Let us consider the $k(\lambda)$-vector space
 $$S(X,Y)_{\lambda}:=\Hom_{k(\lambda)}(X_{\lambda},Y_{\lambda})$$
 provided with the natural structure of $\Gal(K)$-module. Let
 $$A(X,Y,\lambda,K):=\End_{\Gal(K)}(S(X,Y)_{\lambda})$$
 be the centralizer of $\Gal(K)$ in
 $\End_{k(\lambda)}(S(X,Y)_{\lambda})$. Clearly, $A(X,Y,\lambda,K)$ is
 a finite-dimensional $k(\lambda)$-algebra containing the scalars
 $k(\lambda)$.

 \begin{rem}
 \label{basechange2}
Suppose that $L$ is an overfield of $K$ and $K_a$ is the algebraic
closure of $K$ in $L_a$. Let us consider $X$ and $Y$ as abelian
varieties over $L$. It follows from Remark \ref{basechange} that
$\Gal(L)\to \Aut_{k(\lambda)}(S(X,Y)_{\lambda})$ coincides with
the composition  of $\res(L,K):\Gal(L)\to\Gal(K)$ and $\Gal(K)\to
\Aut_{k(\lambda)}(S(X,Y)_{\lambda})$. In particular, the image of
$\Gal(L)\to \Aut_{k(\lambda)}(S(X,Y)_{\lambda})$ lies in the image
of $\Gal(K)\to \Aut_{k(\lambda)}(S(X,Y)_{\lambda})$. It follows
that
$$A(X,Y,\lambda,K)\subset A(X,Y,\lambda,L)\subset
\End_{k(\lambda)}(S(X,Y)_{\lambda}).$$ Clearly, if
$A(X,Y,\lambda,L)$ is a {\sl field} then its every
$k(\lambda)$-subalgebra is also a field, because
$A(X,Y,\lambda,L)$ is finite-dimensional; in particular,
$A(X,Y,\lambda,K)$ is also a field.
 \end{rem}

 \begin{defn}
 \label{disj}
$(X,i)$ and $(Y,j)$ are {\sl disjoint}  at $\lambda$ over $K$if
$A(X,Y,\lambda,K)$ is a field.
 \end{defn}

\begin{rem}
\label{basechange3} It follows from Remark \ref{basechange2} that
if $(X,i)$ and $(Y,j)$ are {\sl disjoint} at $\lambda$ over
$L\supset K$ then they are also disjoint over $K$.
\end{rem}

\begin{thm}
\label{mainA} Suppose that the following conditions hold:
\begin{enumerate}
\item[(i)] The field extensions $K(X_{\lambda})$ and
$K(Y_{\lambda})$  are linearly disjoint over $K$. \item[(ii)]
Consider the centralizer
$k_1:=\End_{\tilde{G}_{\lambda,X}}(X_{\lambda})$ of
$\tilde{G}_{\lambda,X}$ in $\End_{k(\lambda)}(X_{\lambda})$ and
the centralizer $k_2:=\End_{\tilde{G}_{\lambda,Y}}(Y_{\lambda})$
of $\tilde{G}_{\lambda,Y}$ in $\End_{k(\lambda)}(Y_{\lambda})$.
Then the $k(\lambda)$-algebras $k_1$ and $k_2$ are fields that are
linearly disjoint over $k(\lambda)$.
\end{enumerate}
Then $(X,i)$ and $(Y,j)$ are disjoint at $\lambda$ over $K$.
\end{thm}

\begin{thm}
\label{Tdisjoint} If $(X,i)$ and $(Y,j)$ are disjoint at $\lambda$
 over $K$ then one of the following two conditions holds:
\begin{enumerate}
\item[(i)] $\Hom^0((X,i),(Y,j))=0$. \item[(ii)] Both $X$ and $Y$
are isogenous over $K_a$ to  self-products of a certain absolutely
simple abelian variety $B$ of CM-type; in addition, $\End^0(X,i)$
is a $r_X^2$-dimensional central simple $E$-algebra and
$\End^0(Y,i)$ is a $r_Y^2$-dimensional central simple $E$-algebra.
\end{enumerate}
 \end{thm}

We will prove Theorems \ref{mainA} and \ref{Tdisjoint}
  in   \S \ref{proofA}.
We will deduce Theorem \ref{main} from
 the following  statement.

\begin{cor}
\label{mainB} We keep all notations and assumptions of Theorem
\ref{mainA}. Assume in addition that $E$ is normal over $\Q$. Then
one of the following two conditions holds:
\begin{enumerate}
\item[(i)] $\Hom(X,Y)=0, \Hom(Y,X)=0$.
 \item[(ii)] Both $X$ and
$Y$ are isogenous over $K_a$ to  self-products of a certain
absolutely simple abelian variety $B$ of CM-type; in addition,
$\End^0(X,i)$ is a $r_X^2$-dimensional central simple $E$-algebra
and $\End^0(Y,i)$ is a $r_Y^2$-dimensional central simple
$E$-algebra.
\end{enumerate}
\end{cor}

\begin{proof}[of Corollary \ref{mainB}]
Applying Theorems \ref{mainA}  and \ref{Tdisjoint} to
$(X,i\sigma),(Y,j)$ for all $\sigma\in\Gal(E/\Q)$, we conclude
that either the assertion (ii) holds (and we are done)  or all
$$\Hom^0((X,i\sigma),(Y,j))=0.$$ In the latter case, it  follows
from Remark \ref{normal} that $\Hom^0(X,Y)=0$ and therefore
$\Hom(X,Y)=0$, which, in turn, implies that $\Hom(Y,X)=0$.
\end{proof}

\section{Proof of Main Theorem}
\label{proofm} Throughout this section $\ell$ is an odd prime, $K$
a field of characteristic different from $\ell$ and $K_a$ its
algebraic closure,
$$E:=\Q(\zeta_{\ell})\supset \O:=\Z[\zeta_{\ell}]\supset \lambda:=(1-\zeta_{\ell})\cdot
\Z[\zeta_{\ell}], k(\lambda)=\F_{\ell}.$$ Clearly,
$[E:\Q]=\ell-1$.

 Let $f(x) \in K[x]$ be a separable polynomial of degree
$n \ge 4$.

 Let $\RR =\RR_f=\{a_1, \dots , a_n\}\subset K_a$ be the set of
all roots of $f$. We may view the full symmetric group $\Sn$ as
the group of all permutations of  $\RR$. The Galois group
$G=\Gal(f)$ of $f$ permutes the roots and therefore becomes a
subgroup of $\Sn$. The action of $G$ on $\RR$ defines the standard
{{\sl permutational} representation in the $n$-dimensional
$\F_{\ell}$-vector space $\F_p^{\RR}$ of all functions $\psi:\RR
\to \F_{\ell}$. This representation is not irreducible. Indeed,
the "line" of constant functions
 $\F_{\ell}\cdot 1$  and the hyperplane
$(\F_{\ell}^{\RR})^0:=\{\psi\mid \sum_{i=1}^n \psi(a_i)=0\}$ are
$G$-invariant subspaces in $\F_{\ell}^{\RR}$.

 Then we define the {\sl heart}  $(\F_{\ell}^{\RR_f})^{00}$ of the permutational action of
 $G=\Gal(f)$ on
$\RR=\RR_f$ over $\F_{\ell}$ as follows (\cite{Mortimer},
\cite{ZarhinCrelle}).
 If $n$ is not divisible by $\ell$ then we put
   $$(\F_{\ell}^{\RR_f})^{00}=(\F_{\ell}^{\RR})^0:=(\F_{\ell}^{\RR_f})^{0}.$$
If $n$ is divisible by $\ell$ then $(\F_{\ell}^{\RR})^0$ contains
$\F_{\ell}\cdot 1$ and we obtain the natural representation of
$G=\Gal(f)$ in the $(n-2)$-dimensional $\F_{\ell}$-vector
quotient-space $(\F_{\ell}^{\RR})^0/(\F_{\ell}\cdot 1)$. In this
case we put
$$(\F_{\ell}^{\RR_f})^{00}=(\F_{\ell}^{\RR})^{00}:=(\F_{\ell}^{\RR})^0/(\F_{p}\cdot
1).$$ In both cases it is known that the $\Gal(f)$-module
$(\F_{\ell}^{\RR_f})^{00}$ is {\sl faithful} (recall that $n\ge 4$
and $\ell>2$).

\begin{rem}
\label{AS} It is known \cite[Satz 4a]{Klemm} (see also \cite[Lemma
2.4]{ZarhinCrelle}) that if either $n=\deg(f)$ is not divisible by
$\ell$ and $\Gal(f)$ is doubly transitive or $n$ is  divisible by
$\ell$ and $\Gal(f)$ is $3$-transitive  then the centralizer
$\End_{\Gal(f)}((\F_{\ell}^{\RR_f})^{00})=\F_{\ell}$. (Conversely,
one may easily check \cite[Satz 4a]{Klemm} that if $n$ is not
divisible by $\ell$ and $H\subset\Perm(\RR)$ is a permutation
group with $\End_{H}((\F_{\ell}^{\RR)})^{00})=\F_{\ell}$ then $H$
is doubly transitive.)
\end{rem}

\begin{rem}
\label{PS} Let us assume that $K$ contains a primitive $\ell$th
root of unity $\zeta$. Then the map
$$(x,y) \mapsto (x,\zeta y) $$
gives rise to a birational periodic automorphism $\delta_{\ell}$
of $C_{f,\ell}$ with exact period $\ell$. By functoriality,
$\delta_{\ell}$ induces an automorphism of $J(C_{f,\ell})$ which
we still denote by $\delta_{\ell}$. It is known
\cite{Poonen,SPoonen} (see also \cite{ZarhinM}) that
$\delta_{\ell}$ satisfies the $\ell$th cyclotomic equation in
$\End_K(J(C_{f,\ell}))$. This gives rise to the embeddings
$$i_f:\O=\Z[\zeta_{\ell}] \hookrightarrow \End_K(J(C_{f,\ell})), \ E=\Q[\zeta_{\ell}] \hookrightarrow
\End_K^0(J(C_{f,\ell}))$$ with $i_f(1)=1_ {J(C_{f,\ell})}$ and
$i_f(\zeta_{\ell})=\delta_{\ell}$.

Notice that $\lambda=(1-\zeta_{\ell})\cdot \Z[\zeta_{\ell}]$ is
the only maximal ideal dividing $\ell$ in $\Z[\zeta_{\ell}]$ and
the corresponding residue field $k(\lambda)=\F_{\ell}$. The finite
Galois module $J(C_{f,\ell})_{\lambda}$ admits the following
description.
 The
canonical surjection $\Gal(K) \twoheadrightarrow \Gal(f)$
 defines on the $\Gal(f)$-module
$(\F_{\ell}^{\RR_f})^{00}$  the natural structure of
$\Gal(K)$-module. Then the $\F_{\ell}[\Gal(K)]$-modules
$(\F_{\ell}^{\RR_f})^{00}$ and $J(C_{f,\ell})_{\lambda}$ are
canonically isomorphic \cite{Poonen},\cite{SPoonen}. In
particular, this implies that
$$K(J(C_{f,\ell})_{\lambda})=K(\RR_f),$$
(recall that the the $\Gal(f)$-module $(\F_{\ell}^{\RR_f})^{00}$
is {\sl faithful}).
\end{rem}

 Theorem \ref{main} now clearly  is an immediate corollary of
Remark \ref{AS} and the following result.
\begin{thm}
\label{mainC} Suppose that $K$ is a field of characteristic
different from $\ell$ that contains a primitive $\ell$th root of
unity.
  Let $f(x), h(x)\in K[x]$
be separable  polynomials of degree $n\ge 4$ and $m\ge 4$
respectively. Suppose that the splitting fields of
 $f$ and $h$ are linearly disjoint over  $K$.
Suppose that
$$\End_{\Gal(f)}((\F_{\ell}^{\RR_f})^{00})=\F_{\ell}, \ \End_{\Gal(h)}((\F_{\ell}^{\RR_h})^{00})=\F_{\ell}.$$
Then one the two following conditions hold:
\begin{itemize}
\item[(i)]
$\Hom(J(C_{f,\ell}), J(C_{h,\ell}))=0$ and
$\Hom(J(C_{h,\ell}),J(C_{f,\ell}))=0$.
\item[(ii)]
  $p:=\fchar(K)>0$ and
there exists an absolutely simple abelian variety $Z$ defined over
an algebraic closure $\bar{\F}_{p}$ of $\F_{p}$ such that both
$J(C_{f,\ell})$ and $J(C_{h,\ell})$ are abelian varieties of
CM-type isogenous over $K_a$ to self-products of $Z$. In addition,
the centralizer of $\Q[\delta_{\ell}]\cong \Q(\zeta_{\ell})$ in
$\End^0(J(C_{f,\ell}))$ is a central simple
$\Q(\zeta_{\ell})$-algebra of dimension
$\left(\frac{2\dim(J(C_{f,\ell}))}{\ell-1}\right)^2$
 and the centralizer of
$\Q[\delta_{\ell}]\cong \Q(\zeta_{\ell})$ in
$\End^0(J(C_{h,\ell}))$ is a central simple
$\Q(\zeta_{\ell})$-algebra of dimension
$\left(\frac{2\dim(J(C_{h,\ell}))}{\ell-1}\right)^2$.
\end{itemize}
\end{thm}

    \begin{proof}[of Theorem \ref{mainC}]
By the assumption and Remark \ref{PS},
$K(J(C_{f,\ell)})_{\lambda})$ and $K(J(C_{h,\ell)})_{\lambda})$
are linearly disjoint over $K$.

Applying Corollary \ref{mainB} (with $k_1=\F_{\ell}=k_2, \
X=J(C_{f,\ell}), Y=J(C_{h,\ell})$), we conclude that either
$$\Hom(J(C_{f,\ell}),J(C_{h,\ell}))=0, \
\Hom(J(C_{h,\ell}),J(C_{f,\ell}))=0$$ (i.e. the case (i) holds) or
both $J(C_{f,\ell})$ and $J(C_{h,\ell})$ are isogenous over $K_a$
to  self-products of a certain absolutely simple abelian variety
$B$ of CM-type; in addition, the centralizer of $\Q(\zeta_{\ell})$
in $\End^0(J(C_{f,\ell}))$ is a
$\left(\frac{2\dim(J(C_{f,\ell}))}{\ell-1}\right)^2$-dimensional
central simple $\Q(\zeta_{\ell})$-algebra. By \cite[Theorem
3.6]{ZarhinCrelle}, the last property cannot take place  in
characteristic zero and therefore $p:=\fchar(K_a)=\fchar(K)>0$. In
order to check that the case (ii) holds, one has only to recall
that in characteristic $p$ every absolutely simple abelian variety
of CM-type is isogenous to an abelian variety over $\bar{\F}_p$ (a
theorem of Grothendieck \cite{Oort}). By the same token, we get
the desired results for $J(C_{h,\ell})$.
\end{proof}

\begin{rem}
Theorem \ref{mainC} suggests that it may be interesting to
classify subgroups $G=\Gal(f)\subset \Perm(\RR)$ such that
$n=\#(\RR)$ is divisible by $\ell$ and
$\End_{G}((\F_{\ell}^{\RR})^{00})=\F_{\ell}$ (or a field).
According to \cite[Satz 11]{Klemm}, if $G$ is transitive (i.e.
$f(x)$ is irreducible) then such $G$ must be doubly transitive (if
$n\ge 4$ and $\ell$ is odd). The (almost) complete classification
of known doubly transitive $G$ with (absolutely) irreducible
$(\F_{\ell}^{\RR})^{00}$ is given in \cite{Mortimer} (see also
\cite{Ivanov}). Of course, in the irreducible case the centralizer
is a field (and even $\F_{\ell}$ in the absolutely irreducible
case).
\end{rem}

\begin{thm}
\label{mainP} Suppose that $K$ is a field of prime characteristic
different from $\ell$ that contains a primitive $\ell$th root of
unity.
  Let $f(x), h(x)\in K[x]$
be separable  polynomials of degree $n\ge 9$ and $m\ge 4$
respectively. Suppose that the splitting fields of
 $f$ and $h$ are linearly disjoint over  $K$.
Suppose that $\ell$ divides  $n$ and $\Gal(f)$ coincides either
with full symmetric group $\Sn$ or with the alternating group
$\An$. Suppose that
$$ \End_{\Gal(h)}((\F_{\ell}^{\RR_h})^{00})=\F_{\ell}$$
(e.g., $\ell$ does not divide $m$ and $\Gal(h)$ is doubly
transitive.) Then
$$\Hom(J(C_{f,\ell}), J(C_{h,\ell}))=0,\quad \Hom(J(C_{h,\ell}),J(C_{f,\ell}))=0.$$
\end{thm}

\begin{proof}
Clearly, $\Gal(f)$ is $3$-transitive and, thanks to Remark
\ref{AS}, $$\End_{\Gal(f)}((\F_{\ell}^{\RR_f})^{00})=\F_{\ell}.$$
Therefore we may apply Theorem \ref{mainC}. Assume that the
assertion (ii) holds true. In particular, the centralizer of
$\Q[\delta_{\ell}]\cong \Q(\zeta_{\ell})$ in
$\End^0(J(C_{f,\ell}))$ is a central simple
$\left(\frac{2\dim(J(C_{f,\ell}))}{\ell-1}\right)^2$-dimensional
$\Q(\zeta_{\ell})$-algebra. Recall that
$r:=\frac{2\dim(J(C_{f,\ell}))}{\ell-1}=n-1$ or $n-2$; in both
cases we have $r>1$ and therefore the $r^2$-dimensional
centralizer of $\Q[\delta_{\ell}]$ contains an  overfield
$E'\supset \Q[\delta_{\ell}]$ that does {\sl not} coincide with
$\Q[\delta_{\ell}]$. However, it follows from Theorem 0.1 of
\cite{ZarhinSb} that  $\Q[\delta_{\ell}]$ is a maximal commutative
$\Q$-subalgebra in $\End^0(J(C_{f,\ell}))$. This gives us a
desired contradiction.
\end{proof}

\section{Representation theory}
\label{rep}
This Section contains auxiliary results that will be
used in Section \ref{proofA}.
\begin{lem}
\label{irrprod} Let $F$ be a field. Let $H_1$ and $H_2$ be groups.
Let  $\tau_i: H_i \to \Aut_F(W_i)$ $(i=1,2)$ be  linear
finite-dimensional representation of  $H_i$ over $F$ and
$F_i:=\End_{H_i}(W_i)$. Let $W_1^*=\Hom_{F}(W_1,F)$ be the dual of
$W_1$ and $\tau_1^*:H_1\to \Aut_{F}(W_1^*)$ the dual of $\tau_1$.
Let us assume that the $F$-algebras $F_1$ and $F_2$ are fields
that are linearly disjoint over $F$.

Let us consider the natural linear representation
$$\tau_1^*\otimes\tau_2:H_1\times H_2 \to
\Aut_F(\Hom_F(W_1,W_2))$$ of the group $H:=H_1\times H_2$ in the
$F$-vector space $S:=\Hom_F(W_1,W_2)$. Then $\End_H(S)$ is a
field.
\end{lem}

\begin{proof}
 One may easily check that the centralizer of $H_1$ in
$\End_{F}(W_1^*)$ still coincides with $F_1$. Let $A_1$ be the
$F$-subalgebra of $\End_{F}(W_1^*)$ generated by $\tau_1^*(H_1)$;
clearly, the centralizer of $A_1$ in $\End_{F}(W_1^*)$ also
coincides with $F_1$. Similarly, if $A_2$ is the $F$-subalgebra of
$\End_{F}(W_2)$ generated by $\tau_2(H_2)$ then  the centralizer
of $A_2$ in $\End_{F}(W_2)$ coincides with $F_2$. Clearly,   the
$F$-subalgebra of $\End_{F}(W_1^*\otimes_{F}W_2)$ generated by
$\tau_1^*\otimes\tau_2(H_1\times H_2)$ coincides with
$$A_1\otimes_{F}A_2\subset\End_{F}(W_1^*)\otimes_{F}\End_{F}(W_2)=\End_{F}(W_1^*\otimes_{F}W_2).$$
It follows from Lemma (10.37) on p. 252 of \cite{CR} that the
centralizer of $A_1\otimes_{F}A_2$ in
$\End_{F}(W_1^*\otimes_{F}W_2)$ coincides with $F_1\otimes_{F}F_2$
and therefore is a field, thanks to the linear disjointness of
$F_1$ and $F_2$ . This implies that the centralizer of $H_1\times
H_2$ in $\End_{F}(W_1^*\otimes_{F}W_2)$ is the field
$F_1\otimes_{F}F_2$. Since the $H$-modules $W_1^*\otimes_{F}W_2$
and $\Hom_F(W_1,W_2)$ are canonically isomorphic, the centralizer
of $H$ in $\End_{F}(\Hom_F(W_1,W_2))$ is also a field.
\end{proof}

\begin{lem}
\label{irr} Let $L$ be a complete discrete valuation field with
discrete valuation ring $\O_L$, its maximal ideal ${\mathfrak m}$
and residue field $k=\O_L/{\mathfrak m}$. Let $V$ be a
finite-dimensional vector space over $L$, $\tau: G \to \Aut_L(V)$
a completely reducible linear representation of a group $G$ in
$V$. Let $T$ be a $G$-stable $\O_L$-lattice in $V$.  Consider the
finite-dimensional $k$-vector space $\bar{V}=T\otimes_{\O_L}k$
provided with a natural linear representation $\bar{\tau}:G
\to\Aut_k(\bar{V})$ that is the reduction of $\tau$ modulo
${\mathfrak m}$. If the centralizer of $G$ in $\End_k(\bar{V})$ is
a field then $\tau$ is irreducible.
\end{lem}

\begin{proof}
Suppose that $\tau$ is not irreducible. Since it is completely
reducible, there exist non-zero $u_1, u_2 \in \End_G(V)$ with $u_1
u_2=0$. Multiplying (if necessary) both $u_1, u_2$ by suitable
powers of an uniformizer, we may assume that $u_1(T)\subset T,
u_2(T)\subset T$ but neither $u_1$ nor $u_2$ lies in ${\mathfrak
m}\cdot\End_{\O_L}(T)$. It follows that the images $\bar{u}_1,
\bar{u}_2$ of $u_1$ and $u_2$ with respect to the reduction map
$\End_{\O_L}(T)\to \End_k(\bar{V})$ satisfy
$$\bar{u}_1\ne 0, \bar{u}_2\ne 0, \bar{u}_1\bar{u}_2=0.$$
Since both $\bar{u}_1, \bar{u}_2$ obviously lie in the centralizer
of $G$ in $\End_k(\bar{V})$, we get a contradiction.
\end{proof}

\begin{lem}
\label{ext} Let $V$ be a finite-dimensional vector space over a
field $Q$ of characteristic zero, $G$ a group, $\tau:
G\to\Aut_Q(V)$ a completely reducible $Q$-linear representation in
$V$. Let $L$ be an overfield of $Q$ and $i:L \hookrightarrow
\End_Q(V)$ is an embedding of $Q$-algebras that sends $1$ to the
identity automorphism of $V$. Suppose that the image $i(L)$
commutes with $G$. Then the natural $L$-linear representation of
$G$ in $V$ is also completely reducible.
\end{lem}

\begin{proof}
Let $A\subset \End_Q(V)$ be the image of the natural $Q$-algebra
homomorphism $Q[G]\to\End_Q(V)$. The complete reducibility of
$\tau$ means that $A$ is a (finite-dimensional)  semisimple
$Q$-algebra. Therefore $A_L:=A\otimes_Q L$ is a semisimple
$L$-algebra. Clearly, $A\subset \End_L(V)$. This implies that the
image of the natural $L$-algebra homomorphism
$$L[G]\to\End_L(V)\subset \End_Q(V)$$
is isomorphic to a quotient of $A_L$ and therefore is also a
semisimple $L$-algebra. But this means that the natural $L$-linear
representation of $G$ in $V$ is also completely reducible.
\end{proof}

\section{Homomorphisms of abelian varieties: proofs}
\label{proofA}
\begin{proof}[of Theorem \ref{mainA}]
 We need to prove  that the centralizer
 $A(X,Y,\lambda,K)=\End_{\Gal(K)}(S(X,Y)_{\lambda})=\End_{\Gal(K)}(\Hom_{k(\lambda)}(X_{\lambda},Y_{\lambda}))$ of the
 natural representation
$$\Gal(K)
\to\Aut_{k(\lambda)}(\Hom_{k(\lambda)}(X_{\lambda},Y_{\lambda}))$$
is a field. Denote this representation by $\tau$ and let us put
$$F=k(\lambda), H_1=\tilde{G}_{\lambda,X}, W_1=X_{\lambda},
H_2=\tilde{G}_{\lambda,Y}, W_2=Y_{\lambda}.$$ Denote by
$$\tau_1:H_1=\tilde{G}_{\lambda,X}\subset
\Aut_{k(\lambda)}(X_{\lambda})=\Aut_{k(\lambda)}(W_1)$$ and
$$\tau_2:H_2=\tilde{G}_{\lambda,Y}\subset
\Aut_{k(\lambda)}(Y_{\lambda})=\Aut_{k(\lambda)}(W_2)$$ the
corresponding inclusion maps.

It follows from Lemma  \ref{irrprod} that the centralizer of the
linear representation
$$\tau_1^*\otimes\tau_2:\Gal(K(X_{\lambda})/K)\times \Gal(K(Y_{\lambda})/K)
\to\Aut_{k(\lambda)}(\Hom_{k(\lambda)}(X_{\lambda},Y_{\lambda}))$$
 is a field.

One may easily check that   $\tau$, which defines the structure of
$\Gal(K)$-module on $\Hom_{k(\lambda)}(X_{\ell},Y_{\ell})$,
coincides with the composition of the natural surjection  $\Gal(K)
\twoheadrightarrow \Gal(K(X_{\lambda},Y_{\lambda})/K)$, the
natural embedding
$$\Gal(K(X_{\lambda},Y_{\lambda})/K)\hookrightarrow
\Gal(K(X_{\lambda})/K)\times \Gal(K(Y_{\lambda})/K)$$ and
$$\tau_1^*\otimes\tau_2:\Gal(K(X_{\lambda})/K)\times
\Gal(K(Y_{\lambda})/K)\to
\Aut_{k(\lambda)}(\Hom_{k(\lambda)}(X_{\lambda},Y_{\lambda})).$$
Here
 $K(X_{\lambda},Y_{\lambda})$ is the compositum of the fields
 $K(X_{\lambda})$ and $K(Y_{\lambda})$. The linear disjointness of
$K(X_{\lambda})$ and $K(Y_{\lambda})$ means that
$$\Gal(K(X_{\lambda},Y_{\lambda})/K)= \Gal(K(X_{\lambda})/K)\times
\Gal(K(Y_{\lambda})/K).$$ This implies that $\tau$ is the
composition of  the {\sl surjection}  $$\Gal(K) \twoheadrightarrow
\Gal(K(X_{\lambda})/K)\times \Gal(K(Y_{\lambda})/K)$$ and
$\tau_1^*\otimes\tau_2$. Since the centralizer of the
representation
$$\tau_1^*\otimes\tau_2:\Gal(K(Y_{\lambda})/K)\times \Gal(K(X_{\lambda})/K)
\to\Aut_{k(\lambda)}(\Hom_{k(\lambda)}(X_{\lambda},Y_{\lambda}))$$
is a field, the centralizer of the representation $$\tau:
\Gal(K)\to
\Aut_{k(\lambda)}(\Hom_{k(\lambda)}(X_{\lambda},Y_{\lambda}))$$ is
the same field.
\end{proof}

\begin{proof}[of Theorem \ref{Tdisjoint}]
First, we need an additional information about the Tate modules
$T_{\ell}(X)$ and $T_{\ell}(Y)$  of abelian varieties $X$ and $Y$
 \cite{Mumford}. Recall that $T_{\ell}(X)$ and $T_{\ell}(Y)$ are
free $\O_{\lambda}$-modules provided with the continuous actions
of $\Gal(K)$ and one may view may view
$\tilde{\rho}_{\ell,X}:\Gal(K) \to \Aut_{k(\lambda}(X_{\lambda})$
as the reduction of $\rho_{\ell,X}:\Gal(K) \to
\Aut_{\O_{\lambda}}(T_{\ell}(X))$ modulo $\lambda$ and
$\tilde{\rho}_{\lambda,Y}:\Gal(K) \to
\Aut_{k(\lambda}(Y_{\lambda})$ as the reduction of
$\rho_{\ell,Y}:\Gal(K) \to \Aut_{\O_{\lambda}}(T_{\ell}(Y))$
modulo $\lambda$.

It is known \cite{Ribet} that  the Tate $\Q_{\ell}$-modules
$V_{\ell}(X)=T_{\ell}(X)\otimes_{\Z_{\ell}}\Q_{\ell}$ and
$V_{\ell}(Y)=T_{\ell}(Y)\otimes_{\Z_{\ell}}\Q_{\ell}$  are
$\O_{\lambda}\otimes_{\Z_{\ell}}\Q_{\ell}=E_{\lambda}$-vector
spaces of dimension $\frac{2\dim(X)}{[E:\Q]}$ and
$\frac{2\dim(Y)}{[E:\Q]}$ respectively. (Here
$E_{\lambda}=E\otimes_{\Q}\Q_{\ell}$ is the completion of $E$ with
respect to the $\lambda$-adic topology.) The groups $T_{\ell}(X)$
and $T_{\ell}(Y)$ are naturally identified with the
$\O_{\lambda}$-lattices in $V_{\ell}(X)$ and $V_{\ell}(Y)$
respectively and the inclusions
$$\Aut_{\O_{\lambda}}(T_{\ell}(X))\subset
\Aut_{E_{\lambda}}(V_{\ell}(X)), \quad
\Aut_{\O_{\lambda}}(T_{\ell}(Y))\subset
\Aut_{E_{\lambda}}(V_{\ell}(Y))$$ allow us to consider
$V_{\ell}(X)$ and $V_{\ell}(Y)$ as representations of $\Gal(K)$
over $E_{\lambda}$.

Our task now is to prove that the natural representation of
$\Gal(K)$ in
 $$V_1:=\Hom_{E_{\lambda}}(V_{\ell}(Y),V_{\ell}(X))$$
over $E_{\lambda}$ is irreducible. For this, we may and will
assume that $K$ is finitely generated over its prime subfield
(replacing $K$ by a suitable subfield and using Remark
\ref{basechange3}). Then the conjecture of Tate \cite{Tate}
(proven by the author in characteristic $>2$
\cite{ZarhinI,ZarhinT}, Faltings in characteristic zero
\cite{Faltings1,Faltings2} and Mori in characteristic $2$
\cite{MB}) asserts that the natural representation of $\Gal(K)$ in
$V_{\ell}(Z)$ over $\Q_{\ell}$ is completely reducible for any
abelian variety $Z$ over $K$. In particular, the natural
representations of $\Gal(K)$ in $V_{\ell}(X)$ and $V_{\ell}(Y)$
over $\Q_{\ell}$ are completely reducible. It follows from Lemma
\ref{ext}
 that the natural representations of $\Gal(K)$
in $V_{\ell}(X)$ and $V_{\ell}(Y)$ over $E_{\lambda}$ are also
completely reducible.

It follows easily that the dual Galois representation in
$\Hom_{\Q_{\ell}}(V_{\ell}(X),E_{\lambda})$ is also completely
reducible.  Since $E_{\lambda}$ has characteristic zero, it
follows from a theorem of Chevalley \cite[p. 88]{Chevalley} that
the Galois representation in the tensor product
$\Hom_{E_{\lambda}}(V_{\ell}(X),E_{\lambda})\otimes_{E_{\lambda}}V_{\ell}(Y)=
\Hom_{E_{\lambda}}(V_{\ell}(X),V_{\ell}(Y))=V_1$ is completely
reducible.

Second, I claim that the natural representation of  $\Gal(K)$ in
 $V_1$
over $E_{\lambda}$ is irreducible. Indeed, the
$\O_{\lambda}$-module
$\Hom_{\O_{\lambda}}(T_{\ell}(X),T_{\ell}(Y))$ is a
$\Gal(K)$-invariant $\O_{\lambda}$-lattice in
$\Hom_{E_{\lambda}}(V_{\ell}(X),V_{\ell}(Y))=V_1$. On the other
hand, the reduction of this lattice modulo $\lambda$ coincides
with
$$\Hom_{\O_{\lambda}}(T_{\ell}(X),T_{\ell}(Y))\otimes_{\O_{\lambda}}k(\lambda)
= \Hom_{k(\lambda)}(X_{\lambda},Y_{\lambda}).$$ Now the desired
irreducibility follows from Lemma \ref{irr}.

 Third, recall that there is a natural embedding
\cite[Sect. 19]{Mumford}
$$\Hom^0(X,Y)\otimes_{\Q}\Q_{\ell}\subset\Hom_{\Q_{\ell}}(V_{\ell}(X),V_{\ell}(Y)),$$
whose image is a $\Gal(K)$-invariant $\Q_{\ell}$-vector subspace.
Clearly, the image of $\Hom^0((X,i),(Y,j))\otimes_{\Q}\Q_{\ell}$
under this embedding lies in
$\Hom_{E_{\lambda}}(V_{\ell}(X),V_{\ell}(Y))$ and this image is a
$\Gal(K)$-invariant $E\otimes_{\Q}\Q_{\ell}=E_{\lambda}$-vector
subspace of $\Hom_{E_{\lambda}}(V_{\ell}(X),V_{\ell}(Y))$.
 The irreducibility of
$\Hom_{E_{\lambda}}(V_{\ell}(X),V_{\ell}(Y))$ implies that either
$$\Hom^0((X,i),(Y,j))\otimes_{\Q}\Q_{\ell}=\Hom_{E_{\lambda}}(V_{\ell}(X),V_{\ell}(Y))$$
or $\Hom^0((X,i),(Y,j))\otimes_{\Q}\Q_{\ell}=0$. Since
$\Hom^0((X,i),(Y,j))$ is a $E$-vector space, either
$\Hom^0((X,i),(Y,j))=0$ or
$\dim_E(\Hom^0((X,i),(Y,j)))$ equals
$\frac{4\cdot\dim(X)\cdot\dim(Y)}{[E:\Q]^2}$.

In the first case we are done. In the second case the result
follows from Theorem \ref{manyhom}.
\end{proof}


\begin{thebibliography}{99}

\bibitem{Chevalley}  Chevalley C.: Th\'eorie des groupes de Lie,
tome {\bf III}. Hermann,  Paris,  1954.


\bibitem{CR} Curtis, Ch. W., Reiner, I. : Methods of Representation
Theory, Vol. I. John Wiley \& Sons, New York Chichester Brisbane
Toronto, 1981.

\bibitem{DM} Dixon J. D. , Mortimer B. , Permutation Groups. Springer-Verlag,  New York Berlin
Heidelberg, 1996.

\bibitem{Faltings1} Faltings G.: Endlichkeitss\"atze f\"ur
abelsche Variaet\"aten \"uber Z\"ahlkorpern. Invent. Math. {\bf
73}, 349--366 (1983).

\bibitem{Faltings2} Faltings G.: Complements to Mordell. In:
 in: Faltings G. , Wustholz G. et al. (ed.) Rational
points, Chapter VI. Third edition. Aspects of Mathematics, E6.
Friedr. Vieweg \& Sohn, Braunschweig, 1992.

\bibitem{Herstein} Herstein, I. N.: Noncommutative rings.
  Mathematical Association of America, Washington, DC, 1968.

\bibitem{Ivanov} Ivanov, A. A. , Praeger, Ch. E.: On finite
affine 2-Arc transitive graphs. Europ. J. Combinatorics {\bf 14},
421--444 (1993).

\bibitem{Klemm} Klemm, M.: \"Uber die Reduktion von
Permutationsmoduln. Math. Z. {\bf 143}, 113--117 (1975).

\bibitem{MB} Moret-Bailly, L.: Pinceaux de vari\'et\'es ab\'eliennes. Ast\'erisque {\bf 129} (1985).


\bibitem{Mortimer} Mortimer, B.: The modular permutation
representations of the known doubly transitive groups. Proc.
London Math. Soc. (3) {\bf 41}, 1--20 (1980).


\bibitem{Mumford} Mumford, D.: Abelian varieties, Second edition.
 Oxford University Press, London, 1974.

 \bibitem{Oort}Oort,  F.: The isogeny class of a CM-abelian
variety is defined over a finite extension of the prime field. J.
Pure Applied Algebra {\bf 3}, 399--408 (1973).



\bibitem{Poonen}Poonen,  B., Schaefer, E.: Explicit descent for Jacobians of cyclic covers of
 the projective line. J. reine angew. Math. {\bf 488}, 141--188 (1997).

 \bibitem{Ribet}Ribet,  K.: Galois action on division points of Abelian varieties with real
multiplications. Amer. J. Math. {\bf 98}, 751--804 (1976).


\bibitem{SPoonen} Schaefer, E.: Computing a Selmer group of a Jacobian using functions on the curve.
 Math. Ann. {\bf 310}, 447--471 (1998).

\bibitem{Schur} Schur, I.: Gleichungen ohne Affect.
Sitz. Preuss. Akad. Wiss., Physik-Math. Klasse 443--449 (1930) (=
Ges. Abh. III, 191--197).

\bibitem{SerreM} Serre, J.-P.: Lectures on the Mordell-Weil
Theorem, 2nd edition, Friedr. Vieweg \& Sons,
Braunschweig/Wiesbaden, 1990.

\bibitem{SerreG} Serre, J.-P.: Topics in Galois Theory,
Jones and Bartlett Publishers, Boston-London, 1992.

\bibitem{Serre} Serre, J.-P.: Repr\'esentations lin\'eares des groupes
finis, Troisi\'eme \'edition. Hermann, Paris, 1978.

\bibitem{Tate} Tate, J. Endomorphisms of Abelian varieties over finite
fields. Invent. Math. {\bf 2}, 134--144 (1966).


\bibitem{ZarhinI} Zarhin, Yu. G.: Endomorphisms of Abelian varieties over fields of finite
characteristic. Izv. Akad. Nauk SSSR ser. matem. {\bf 39},
272--277 (1975); Math. USSR Izv. {\bf 9}, 255 - 260 (1975).

\bibitem{ZarhinT} Zarhin, Yu. G.: Abelian varieties in characteristic
$P$. Mat. Zametki {\bf 19}, 393--400 (1976); Mathematical Notes
{\bf 19}, 240--244 (1976).


\bibitem{ZarhinCrelle} Zarhin, Yu. G.: Cyclic covers, their Jacobians and endomorphisms.  J.
reine angew.  Math. {\bf 544}, 91--110 (2002).

\bibitem{ZarhinCamb} Zarhin, Yu. G.: The endomorphism rings of Jacobians of cyclic covers of the
projective line. Math. Proc. Cambridge Philos. Soc. {\bf 136},
257--267 (2004).

\bibitem{ZarhinSb} Zarhin, Yu. G.: Endomorphism rings of certain Jacobians in finite
characteristic. Matem. Sbornik  {\bf 193}, issue 8, 39--48 (2002);
Sbornik Math. {\bf 193} (8), 1139-1149 (2002).

\bibitem{ZarhinSh} Zarhin,  Yu. G.:  Homomorphisms of
hyperelliptic Jacobians.  In: Number Theory, Algebra and Algebraic
Geometry (Shafarevich Festschrift). Trudy Mat. Inst. Steklov {\bf
241}, 90--104 (2003); Proc. Steklov Inst. Math. {\bf 241}, 79--92
(2003).

\bibitem{ZarhinL} Zarhin, Yu. G.:  Homomorphisms of abelian varieties.
In: Y. Aubry, G. Lachaud (ed.)
 Arithmetic, Geometry and Coding Theory (AGCT 2003),
 S\'eminaires et Congr\'es {\bf 11}, 189--215 (2005).


\bibitem{ZarhinM} Zarhin,  Yu. G.: Endomorphism algebras of superelliptic Jacobians. In:
 F. Bogomolov, Yu. Tschinkel (ed.) Geometric methods in Algebra and Number Theory,
Progress in Math. {\bf 235}, 339--362, Birkh\"auser, Boston Basel
Berlin, 2005.
\end{thebibliography}
\end{document}